\documentclass[10pt]{IEEEtran}

\usepackage{xspace,amssymb,epsfig,subfigure,syntonly}
\usepackage{epsfig,amsmath,color}
\usepackage{xspace,syntonly}

\newcommand{\utwi}[1]{\mbox{\boldmath $#1$}}

\newcommand{\trace}{{\textrm{Tr}}}
\newcommand{\rank}{{\textrm{rank}}}

\newcommand{\cD}{{\cal D}}

\newcommand{\cN}{{\cal N}}
\newcommand{\cP}{{\cal P}}

\newcommand{\cS}{{\cal S}}

\newcommand{\cT}{{\cal T}}

\newcommand{\cE}{{\cal E}}
\newcommand{\cI}{{\cal I}}

\newcommand{\cH}{{\cal H}}
\newcommand{\cV}{{\cal V}}

\newcommand{\ba}{{\bf a}}

\newcommand{\be}{{\bf e}}

\newcommand{\bx}{{\bf x}}
\newcommand{\bu}{{\bf u}}
\newcommand{\bv}{{\bf v}}
\newcommand{\bi}{{\bf i}}

\newcommand{\bB}{{\bf B}}

\newcommand{\bD}{{\bf D}}

\newcommand{\bT}{{\bf T}}

\newcommand{\bX}{{\bf X}}
\newcommand{\bZ}{{\bf Z}}
\newcommand{\bY}{{\bf Y}}
\newcommand{\bV}{{\bf V}}

\newcommand{\bPhi}{{\utwi{\Phi}}}

\begin{document}

\newtheorem{definition}{Definition}
\newtheorem{remark}{Remark}
\newtheorem{proposition}{Proposition}
\newtheorem{lemma}{Lemma}

%------------------------------------------------------------------------------
% Title.
%------------------------------------------------------------------------------
\title{Economic Dispatch in Unbalanced Distribution Networks via Semidefinite Relaxation}

\author{Emiliano Dall'Anese, Georgios B. Giannakis, and Bruce F. Wollenberg
\thanks{\protect\rule{0pt}{0.5cm}%
Version of the manuscript: June 29, 2012. 
The authors are with the Department of Electrical and Computer Engineering, University of Minnesota, 200 Union Street SE, %
Minneapolis, MN 55455, USA. 
E-mails: {\tt \{emiliano, georgios, wollenbe\}@umn.edu}.
Part of this work was presented to the 44th North
America Power Symposium, Urbana-Champaign, IL, USA.
}
}

% The paper headers
\markboth{}%
{First Author \MakeLowercase{\textit{et al.}}: Title}

% make the title area
\maketitle

%\thispagestyle{empty}

%%%%%%%%%%%%%%%%%%%%%%%%%%%%%%%%%%%%%%%%%%%%%%%
% Abstract
%%%%%%%%%%%%%%%%%%%%%%%%%%%%%%%%%%%%%%%%%%%%%%%
\begin{abstract}
The economic dispatch problem is considered for unbalanced three-phase power distribution networks entailing both non-deferrable and elastic loads, and distributed generation (DG) units. The objective is to minimize the costs of power drawn from the main grid and supplied by the DG units over a given time horizon, while meeting the overall load demand and effecting voltage regulation. Similar to optimal power flow counterparts for balanced systems, the resultant optimization problem is nonconvex. Nevertheless, a semidefinite programming (SDP) relaxation technique is advocated to obtain a (relaxed) convex problem solvable in polynomial-time complexity. To promote a reliable yet efficient feeder operation, SDP-compliant constraints on line and neutral current magnitudes are accommodated in the optimization formulated, along with constraints on the power factor at the substation and at nodes equipped with capacitor banks. Tests on the IEEE 13-node radial feeder demonstrate the ability of the proposed method to attain the \emph{globally} optimal solution of the original nonconvex problem.
\end{abstract}


%%%%%%%%%%%%%%%%%%%%%%%%%%%%%%%%%%%%%%%%%%%%%%%
\section{Introduction}
\label{sec:Introduction}
%%%%%%%%%%%%%%%%%%%%%%%%%%%%%%%%%%%%%%%%%%%%%%%

The advent of distributed energy resources, along with the rapid proliferation of controllable loads such as, e.g., plug-in hybrid electric vehicles (PHEVs), call for innovative energy management methodologies to ensure highly efficient operation of distribution networks, effect voltage regulation, and facilitate emergency response~\cite{Momoh09}. Toward these goals, variants of the optimal power flow (OPF) problem have been devised with the objective of optimizing the power supplied by distributed generation (DG) units as well as by the utility at the substation, subject to electrical network constraints on powers and voltages, and the expected load profile~\cite{Khodr07}.  

These approaches however, are deemed challenging because they require solving nonconvex problems. Non-convexity stems from the nonlinear relationship between voltages and the apparent powers demanded at the loads. Furthermore, the high resistance-reactance ratio in conventional distribution lines severely challenges the convergence of Newton-Raphson iterations, which have been traditionally employed for solving nonconvex OPF problems in transmission networks~\cite{Irving87,Tripathy82}. This has motivated the adoption of forward/backward sweeping methods~\cite{Cespedes90}, which enable computationally-efficient load flow analysis, but are not suitable for optimization purposes, fuzzy dynamic programming~\cite{Lu97}, particle swarm optimization~\cite{Sortomme09},  sequential quadratic optimization~\cite{Driesen10}, and steepest descent-based methods~\cite{Forner}. However, these approaches generally return sub-optimal load flow solutions, and may be computationally cumbersome. To alleviate these concerns, the semidefinite programming (SDP) reformulation of~\cite{Bai08} and \cite{LavaeiLow} was recently extended to \emph{balanced} distribution networks in~\cite{Tse12}, and conditions ensuring global optimality of the obtained solution were derived. 

Three-phase distribution feeders however, are inherently \emph{unbalanced} because \emph{i)} a large number of unequal single-phase loads must be served, and \emph{ii)} non-equilateral conductor spacings
of three-phase line segments are involved~\cite{Kerstingbook}. As a consequence, optimization approaches can not rely on single-phase equivalent models (as in, e.g.~\cite{Tse12, Forner}). For the unbalanced setup, an OPF framework was proposed in~\cite{Paudyal11}, where commercial solvers of nonlinear programs were used, and in~\cite{Bruno11}, where Newton methods where utilized in conjunction with OpenDSS load flow solvers. A model based on sequence components was adopted in~\cite{Rashid05}, and the Newton-Raphson algorithm was used. 
However, since these methods are inherently related to gradient descent solvers of nonconvex programs, they inherit the limitations of being sensitive to initialization, and do not guarantee global optimality of their solutions.

The main contribution of the present paper consists in permeating the benefits of SDP relaxation techniques to the economic dispatch problem for \emph{unbalanced} three-phase power distribution systems. This powerful optimization tool not only offers the potential of finding the \emph{globally} optimal solution of the original nonconvex problem in \emph{polynomial-time} complexity~\cite{luospmag10}, but also facilitates the introduction of thermal and quality-of-power constraints without exacerbating the problem complexity. The focus here is on the case where the costs of power provided by the
utility company and supplied by the DG units are known in advance over a given time horizon. Then, the goal is to minimize the overall energy cost so that both non-deferrable and elastic load demands are met, and the node voltages stay within prescribed limits. Furthermore, constraints on line and neutral current magnitudes, as well as on the power factor at substation and nodes equipped with capacitor banks are accommodated in the optimization problem in order to improve reliability and efficiency of the distribution feeder.

{{\it Notation:} Upper (lower) boldface
letters will be used for matrices (column vectors); $(\cdot)^\cT$ for transposition; $(\cdot)^*$ complex-conjugate; and, $(\cdot)^\cH$ complex-conjugate transposition;  
$\Re\{\cdot\}$ denotes the real part, and $\Im\{\cdot\}$ the imaginary part; $j = \sqrt{-1}$ represents the imaginary unit. $\trace(\cdot)$ denotes  the matrix trace; $\rank(\cdot)$ the matrix rank; $\circ$ the Hadamard product; and, $|\cdot|$ denotes the magnitude of a number or the cardinality of a set. 
Given a vector $\bv$ and matrix $\bV$, $[\bv]_{\cP}$ denotes a $|\cP| \times 1$ sub-vector containing the entries of $\bv$ indexed by the set $\cP$, and $[\bV]_{\cP_1,\cP_2}$ the $|\cP_1| \times |\cP_2|$ sub-matrix with row and column indexes described by $\cP_1$ and $\cP_2$. Finally, $\mathbf{0}_{M\times N}$ and $\mathbf{1}_{M\times N}$ denotes $M \times N$ matrices with all zeroes and ones, respectively.

%%%%%%%%%%%%%%%%%%%%%%%%%%%%%%%%%%%%%%%%%%%%%%%
\section{Modeling and problem formulation}
\label{sec:Modeling}
%%%%%%%%%%%%%%%%%%%%%%%%%%%%%%%%%%%%%%%%%%%%%%%

Consider a radial distribution feeder comprising $N$ nodes collected
in the set $\cN := \{1,\ldots,N\}$, and overhead or underground lines represented by the set of edges $\cE := \{(m,n)\} \subset (\cN \cup \{0\})  \times (\cN \cup \{0\})$,  
where the additional node $0$ corresponds to the point of common coupling (PCC). 
The feeder operation is to be optimized over a given time interval $\cI := \{1,2,\ldots,T\}$, where each time slot can represent 
e.g., ten minutes or fifteen  minutes, one hour, etc, depending on the specific short-, medium-, or long-range scheduling horizon~\cite{Paudyal11}. 

The backbone of the feeder generally consists of three-phase lines, with two- and single-phase connections at times present on laterals and sub-laterals. Let $\cP_{mn} \subseteq \{a,b,c\}$ and $\cP_{n} \subseteq \{a,b,c\}$ denote the phase of line $(m,n) \in \cE$ and node $n \in \cN$, respectively. Further, let $V_{n,t}^{\phi} \in \mathbb{C}$ and $I_{n,t}^{\phi} \in \mathbb{C}$ be the complex line-to-ground voltage at node $n \in \cN$ and time slot $t$ of phase $\phi \in \cP_n$, and the current injected at the same node, phase, and time. As usual, the voltages $\bv_{0,t} :=  [V_{0,t}^{a},V_{0,t}^{b},V_{0,t}^{c}]^{\cT}$ at the PCC are assumed to be known~\cite{Kerstingbook}. 

Lines $(m,n) \in \cE$ are modeled as $\pi$-equivalent components, and the $|\cP_{mn}| \times |\cP_{mn}|$ phase impedance and shunt admittance matrices are denoted as $\bZ_{mn} \in \mathbb{C}^{|\cP_{mn}| \times |\cP_{mn}|}$ and $\bY_{mn}^{(s)} \in \mathbb{C}^{|\cP_{mn}| \times |\cP_{mn}|}$, respectively. %\footnote{Impedances and admittances are generally collected in $3\times3$ matrices; in case of two- or one-phase lines, $|\cP_{mn}| \times |\cP_{mn}|$ sub-matrices are extracted.}.
%Notice that, differently from balanced transmission and distribution systems, these matrices are not in general scaled versions of the identity matrix, and they capture self and mutual inductive reactances and capacitive coefficients, respectively, in their off-diagonal elements~\cite[Ch.~6]{Kerstingbook}. 
If four-wire grounded wye lines or lines with multi-grounded neutrals are present, matrices $\bZ_{mn}$ and $\bY_{mn}^{(s)}$ can be obtained from the higher-dimensional ``primitive'' matrices via Kron reduction~\cite[Ch. 6]{Kerstingbook}. Using the $\pi$-equivalent model, it follows from Kirchhoff's current law that the current $I_n^{\phi}$ can be expressed as~\cite{Kerstingbook} 
%\begin{align}
%I_{n,t}^{\phi} & = \sum_{m \in \cN_m}  \left[ \left(\frac{1}{2}\bY_{mn}^{(s)} + \bZ_{mn}^{-1} \right) [\bv_{n,t}]_{\cP_{mn}} \right. \nonumber \\
%& \hspace{3cm} - \bZ_{mn}^{-1} [\bv_{m,t}]_{\cP_{mn}} \Big]_{\{\phi\}}
%\label{currents}
%\end{align}

\vspace{-.3cm}
\small
\begin{align}
I_{n,t}^{\phi} & = \sum_{m \in \cN_m}  \left[ \left(\frac{1}{2}\bY_{mn}^{(s)} + \bZ_{mn}^{-1} \right) [\bv_{n,t}]_{\cP_{mn}}  - \bZ_{mn}^{-1} [\bv_{m,t}]_{\cP_{mn}} \right]_{\{\phi\}}
\label{currents}
\end{align}
\normalsize
where $\cN_n \subset \cN$ is the set of nodes linked to $n$ through a transmission line, and $\bv_{n,t} \in \mathbb{C}^{|\cP_{n}|}$ denotes the column vector collecting the voltages at node $n$ and time slot $t$. Three- or single-phase transformers (if any) are modeled as series components with transmission parameters that depend on the connection type~\cite[Ch.~8]{Kerstingbook},~\cite{Paudyal11}. 

%%%%%%%%%%%%%%%%%%%%
%\begin{figure}[t]
%\begin{center}
%\includegraphics[width=0.5\textwidth]{F_3Phase}
%\caption{Feeder section with 3- and 2-phase lines, and its $\pi$-equivalent model.}
%\label{fig:linemodel}
%\end{center}
%\vspace{-.5cm}
%\end{figure}
%%%%%%%%%%%%%%%%%%%%
Per phase $\phi \in \cP_n$ and node $n \in \cN$, the following two classes of loads are considered. 
\begin{itemize}
	\item A base \emph{non-deferrable} load with active and reactive powers demanded at time $t$ denoted by $P_{L,n,t}^{\phi}$ and $Q_{L,n,t}^{\phi}$, respectively. 
	\item A set $\cD_n^{\phi}$ of \emph{controllable} (elastic) loads, each with a prescribed energy requirement $E_{d,n}^{\phi}$  to be completed over a given interval $\cI_{d,n}^{\phi} := [s_{d,n}^{\phi},f_{d,n}^{\phi}] \subseteq \cI$, with $s_{d,n}^{\phi}$ representing the starting time, and $f_{d,n}^{\phi}$ the termination slot; that is,
%\begin{align}
$\sum_{t \in \cI_{d,n}^{\phi}} \bar{P}_{d,n,t}^{\phi} \Delta_t = E_{d,n}^{\phi}$, 
%\label{controllableload}
%\end{align}
with $\bar{P}_{d,n,t}^{\phi}$ the amount of active power supplied to the controllable load $d \in \cD_n^{\phi}$ at time slot $t$, and $\Delta_t > 0$ the duration of the time slot. 
\end{itemize}
An example of controllable load is PHEVs, whose charging process can be shifted from hours with high price of electricity~\cite{Aliprantis12}, and high load conditions of the distribution network~\cite{Driesen10, Sundstrom12}, to off-peak hours. In this case, users specify the time when PHEVs will be plugged in, and the time by which the charging has to be completed~\cite{Driesen10, Aliprantis12}. 
%Loads are assumed wye-connected; in case of delta-connected load, a delta-to-wye transformation can be employed as usual.  

In distribution feeders, capacitor banks are mounted at some selected nodes to provide reactive power support, aid voltage regulation, and correct the load power factor (PF). As usual, capacitors can be modeled as wye or delta loads with constant susceptance~\cite{Paudyal11},~\cite[Ch.~9]{Kerstingbook}. Therefore, with $y_{C,n}^{\phi}$ denoting the susceptance of a capacitor connected at node $n$ and phase $\phi$,
the reactive power $Q_{C,n,t}^{\phi}$ provided by the capacitor at time $t$ is given by 
%\begin{equation}
$Q_{C,n,t}^{\phi} = y_{C,n}^{\phi} |V_{n,t}^{\phi}|^2$. 
%\end{equation}
To satisfy the load demand, DG units such as, e.g. diesel generators and fuel cells can be employed to complement the power drawn from the main distribution grid. Then, suppose that $S$ DG units are located at nodes $\cS \subset \cN$, and let $P_{G,s,t}^{\phi}$ and $Q_{G,s,t}^{\phi}$ denote the active and reactive powers supplied by unit $s \in  \cS$.  
  
The focus is on the case where the costs of power provided by the
utility company $P_{0,t}^{\phi} := \Re\{V_{0,t}^{\phi} (I_{0,t}^{\phi})^*\}$, $\phi = a,b,c$, 
and supplied or consumed by the DG units are determined in advance for the period $\cI$. Let $\{\kappa_{0,t}\}$ denote the former, and $\{c_{s,t}\}$ the latter.
Then, the goal is to minimize the overall cost of power purchased from the main grid and generated within the feeder (economic dispatch),  so that the total load demand is met, and the node voltages stay within prescribed limits; that is, the following problem is to be solved, where $\cV^{(1)} := \{\{I_{0,t}^{\phi}\}_{\forall \phi,t}, \{V_{n,t}^{\phi}, I_{n,t}^{\phi}, P_{G,n,t}^{\phi},Q_{G,n,t}^{\phi}\}_{\forall \, \phi, n, t}\}$ and $\cV_d :=  \{ \{\bar{P}_{d,n,t}^{\phi} \}_{\forall \phi,n,t}\} $ collect the optimization variables:
\begin{subequations}
\label{Pmg}
\begin{align} 
& \hspace{-.5cm} (P1)\,\, \min_{\substack{\cV^{(1)}\\ \cV_d }} \sum_{t \in \cI} \left(\kappa_{0,t} \sum_{\phi \in \cP_0} P_{0,t}^{\phi} + \sum_{s \in \cS}  c_{s,t} \sum_{\phi \in \cP_s} P_{G,s,t}^{\phi} \right) \label{mg-cost} \\
\textrm{s.t.} \,\, 
& V_{s,t}^{\phi} (I_{s,t}^{\phi})^* = P_{G,s,t}^{\phi} - P_{L,s,t}^{\phi} - \sum_{d \in \cD_s^{\phi}} \bar{P}_{d,s,t}^{\phi} \nonumber \\
&  \hspace{.7cm}+ j Q_{G,s,t}^{\phi} -  j Q_{L,n,t}^{\phi}, \,\,    \forall \,\,  t \in \cI, \phi \in \cP_s, \, s \in \cS \label{mg-balance-source} \\
& V_{n,t}^{\phi} (I_{n,t}^{\phi})^* = - P_{L,n,t}^{\phi} - \sum_{d \in \cD_n^{\phi}} \bar{P}_{d,n,t}^{\phi} - j Q_{L,n,t}^{\phi} \nonumber \\ 
& \hspace{.7cm}  + j y_{C,n}^{\phi} |V_{n,t}^{\phi}|^2, \hspace{.2cm} \forall\,\,  t \in \cI, \phi \in \cP_n, \,  n \in \cN \backslash \cS \label{mg-balance} \\
& \sum_{t \in \cI_{d,n}^{\phi}} \bar{P}_{d,n,t}^{\phi} \Delta_t = E_{d,n}^{\phi}, \hspace{.1cm}  \forall \,\, d \in \cD_n^{\phi}, \phi \in \cP_n, n \in \cN  \label{mg-controllable} \\
&  \hspace{.6cm} 0 \leq \bar{P}_{d,n,t}^{\phi} \leq  \bar{P}_{d,n}^{\textrm{max}}, \hspace{.2cm} \forall \,\,  t \in \cI, d \in \cD_n^{\phi},\,  n \in \cN \label{mg-Pdlimits} \\
& V_{n,t}^{\mathrm{min}} \leq |V_{n,t}^{\phi}| \leq V_{n,t}^{\mathrm{max}} , \hspace{.15cm} \forall \,\,  t \in \cI, \phi \in \cP_n,\,  n \in \cN  \label{mg-Vlimits} \\
& P_{G,s}^{\textrm{min}} \leq P_{G,s,t}^{\phi} \leq  P_{G,s}^{\textrm{max}}, \hspace{.25cm} \forall \,\,  t \in \cI, \phi \in \cP_s,\,  s \in \cS \label{mg-plimits} \\
& Q_{G,s}^{\textrm{min}} \leq Q_{G,s,t}^{\phi} \leq  Q_{G,s}^{\textrm{max}}, \hspace{.15cm} \forall \,\,  t \in \cI, \phi \in \cP_n,\,  n \in \cN \label{mg-qlimits} 
\end{align}
\end{subequations}
where $P_{G,s}^{\textrm{min}},P_{G,s}^{\textrm{max}}, Q_{G,s}^{\textrm{min}},Q_{G,s}^{\textrm{max}}$ capture physical and operational constraints of the DG units, and $V_n^{\mathrm{min}}$ and $V_n^{\mathrm{max}}$ are given minimum and maximum utilization and service voltages. Finally, $\bar{P}_{d,n}^{\textrm{max}}$ represents a possible cap for $\bar{P}_{d,n}^{\phi}$.  If capacitor banks are not present,~\eqref{mg-balance} should be modified to $V_{n,t}^{\phi} (I_{n,t}^{\phi})^* = - P_{L,n,t}^{\phi} - jQ_{L,n,t}^{\phi}$. Recall that the voltages at the PCC are assumed known. However, if needed, this assumption can be relaxed, and (P1) can be appropriately re-stated.

Unfortunately, (P1) is a nonlinear \emph{nonconvex} problem due to the load flow equations~\eqref{mg-balance-source}-\eqref{mg-balance} as well as the voltage constraints~\eqref{mg-Vlimits}. In the next section, an equivalent reformulation of (P1) will be derived, and its solution will be tackled by employing an SDP relaxation technique. %The merit of the proposed solution approach consists in offering the potential of finding the \emph{globally optimal} solution of (P1) in \emph{polynomial-time} complexity~\cite{luospmag10}.  

%%%%%%%%%%%%%%%%%%%%%%%%%%%%%%%%%%%%%%%%%%%%%%%
\section{Relaxed semi-definite programming}
\label{sec:semidefinite}
%%%%%%%%%%%%%%%%%%%%%%%%%%%%%%%%%%%%%%%%%%%%%%%

Consider first a distribution feeder with only three-phase lines and nodes; that is, $|\cP_{n}| = 3$ for all $n \in \cN$, and $|\cP_{nl}| = 3$ for all lines $(n,l) \in \cE$. Let $\mathbf{Y} \in \mathbb{C}^{3(N+1) \times 3(N+1)}$ be a symmetric matrix defined as [cf.~\eqref{currents}]
\begin{equation}
[\bY]_{\cP_n,\cP_m} := \left \{\begin{array}{ll} \sum_{m \in \cN_m}  \left(\frac{1}{2}\bY_{mn}^{(s)} + \bZ_{mn}^{-1} \right), & \textrm{if } m=n\\ 
-\bZ_{mn}^{-1}, &\textrm{if }(m,n)\in \cE\\
\mathbf{0}_{3 \times 3}, & \textrm{otherwise} \end{array} \right.
\label{eq:Ymatrix} \nonumber
\end{equation}
and define the $3(N+1) \times 1$ vectors $\bv_t := [\bv_{0,t}^\cT, \ldots, \bv_{N,t}^\cT]^\cT$ and $\bi_t := [\bi_{0,t}^\cT, \ldots, \bi_{N,t}^\cT]^\cT$, with $\bi_{n,t} := [I^a_{n,t},I^b_{n,t},I^c_{n,t}]^\cT$.
Then,~\eqref{currents} can be re-written in vector-matrix form as $\bi_{t} = \bY \bv_{t}$.  

Since the PCC voltages $\{\bv_{0,t}\}$ are known, re-write the vector of complex voltages as $\bv_{t} = \ba_{0,t} \circ \bx_{t}$, with $\bx_{t}:= [\mathbf{1}_3^\cT, \bv_{1,t}^\cT, \ldots, \bv_{N,t}^\cT]^\cT$ and $\ba_{0,t} := [\bv_{0,t}^\cT,\mathbf{1}_{|\cP_{1}|}^{\cT},\ldots,\mathbf{1}_{|\cP_{N}|}^{\cT}]^\cT$, for all $t \in \cI$. Then, consider expressing the active and reactive powers injected at each node at time $t$, as well as the voltage magnitudes as linear functions of the outer-product matrix $\bX_t := \bx_t \bx_t^\cH$. To this end, define the following admittance-related matrix per node $n$ and phase $\phi$
\begin{equation}
\bY_n^{\phi} := \bar{\be}_n^{\phi} (\bar{\be}_n^{\phi})^{\cT} \bY 
\label{eq:Ynode}
\end{equation}
where $\bar{\be}_n^{\phi} := [\mathbf{0}_{|\cP_0|}^{\cT},\ldots,\mathbf{0}_{|\cP_{n-1}|}^{\cT},\be^{\phi, \cT},\mathbf{0}_{|\cP_{n+1}|}^{\cT},\ldots,\mathbf{0}_{|\cP_N|}^{\cT}]^{\cT}$, and $\{\be^\phi\}_{\phi \in \{a,b,c\}}$ denotes the canonical basis of $\mathbb{R}^3$.  Denote for future use the Hermitian matrices   
\begin{subequations}
\label{eq:Phi}
\begin{align}
\bPhi_{P,n,t}^{\phi} & := \frac{1}{2}\bD_{0,t}^{\cH} (\bY_n^{\phi} + (\bY_n^{\phi})^\cH) \bD_{0,t} \label{eq:PhiP} \\
\bPhi_{Q,n,t}^{\phi} & := \frac{j}{2}\bD_{0,t}^{\cH} (\bY_n^{\phi} - (\bY_n^{\phi})^\cH) \bD_{0,t} \label{eq:PhiQ} \\
\bPhi_{V,n,t}^{\phi} & := \bD_{0,t}^{\cH} \bar{\be}_n^{\phi} (\bar{\be}_n^{\phi})^{\cT} \bD_{0,t}  \label{eq:PhiV} 
\end{align}
\end{subequations}
with $\bD_{0,t} := \textrm{diag}(\ba_{0,t})$. Using~\eqref{eq:Phi}, a linear model in $\bX_t$ (and therefore in $\bV_t := \bv_t \bv_t^\cH$) is established in the following lemma (see also~\cite{ZhuNAPS11} and~\cite{LavaeiLow}).

%\vspace{.2cm}

\begin{lemma}
\label{SDPreformulation}
Apparent powers and voltage magnitudes are linearly related with $\{\bX_t\}$ as [cf.~\eqref{eq:Phi}] 
\begin{subequations}
 \label{eq:PQV} 
\begin{align}
\trace(\bPhi_{P,n,t}^{\phi} \bX_t) & = P_{G,n,t}^{\phi} - P_{L,n,t}^{\phi} - \sum_{d \in \cD_s^{\phi}} \bar{P}_{d,s,t}^{\phi} \label{eq:P} \\
\trace(\bPhi_{Q,n,t}^{\phi} \bX_t) & = Q_{G,n,t}^{\phi} - Q_{L,n,t}^{\phi} + y_{C,n}^{\phi}\trace(\bPhi_{V,n,t}^{\phi} \bX_t)  \label{eq:Q}  \\ 
\trace(\bPhi_{V,n,t}^{\phi} \bX_t) & = |V_{n,t}^{\phi}|^2  \label{eq:V} 
\end{align}
\end{subequations}
with $P_{G,n,t}^{\phi} = Q_{G,n,t}^{\phi} = 0$ for $n \in \cN\backslash \cS$, and $y_{C,n}^{\phi} = 0$ if capacitor banks are not present at node $n$. 
\end{lemma}
\emph{Proof.}   See the Appendix.  \hfill $\Box$

%\vspace{.2cm}

Using~\eqref{eq:P}--\eqref{eq:V}, problem (P1) is \emph{equivalently} reformulated as follows: 
\begin{subequations}
\label{mg2}
\begin{align} 
& \hspace{-.5cm} (P2) \,\, \min_{\{\bX_t\},\cV_d} \sum_{t \in \cI} \kappa_{0,t} \sum_{\phi \in \cP_0} \trace(\bPhi_{P,0,t}^{\phi} \bX_t) \nonumber \\ 
& \hspace{1.1cm} + \sum_{t \in \cI}  \sum_{s \in \cS}  c_{s,t} \sum_{\phi \in \cP_s} \trace(\bPhi_{P,s,t}^{\phi} \bX_t) \\
\textrm{s.t.} \,\, 
& \trace(\bPhi_{P,n,t}^{\phi} \bX_t) + P_{L,n,t}^{\phi} + \sum_{d \in \cD_s^{\phi}} \bar{P}_{d,s,t}^{\phi} = 0, \nonumber \\ 
&\hspace{3cm} \forall\,\, t \in \cI, \phi \in \cP_n, \, \forall \, n \in \cN \backslash \cS  \label{mg2-P} \\
& \trace(\bPhi_{Q,n,t}^{\phi} \bX_t) + Q_{L,n,t}^{\phi} -y_{C,n}^{\phi}\trace(\bPhi_{V,n}^{\phi} \bX_t)  = 0, \nonumber \\
& \hspace{3cm} \forall\,\, t \in \cI, \phi \in \cP_n, \, \forall \, n \in \cN \backslash \cS \label{mg2-Q} \\
&P_{G,s}^{\textrm{min}} \leq \trace(\bPhi_{P,s,t}^{\phi} \bX_t) + P_{L,s,t}^{\phi} + \sum_{d \in \cD_s^{\phi}} \bar{P}_{d,s,t}^{\phi}\leq P_{G,s}^{\textrm{max}},  \nonumber \\
& \hspace{3cm} \forall\,\, t \in \cI, \phi \in \cP_s, \, \forall \, s \in \cS   \label{mg2-Pg} \\
&Q_{G,s}^{\textrm{min}} \leq \trace(\bPhi_{Q,s,t}^{\phi} \bX_t) + Q_{L,s,t}^{\phi} \leq Q_{G,s}^{\textrm{max}},  \nonumber \\
& \hspace{3cm} \forall\,\, t \in \cI, \phi \in \cP_s, \, \forall \, s \in \cS   \label{mg2-Pg} \\
& (V_n^{\mathrm{min}})^2 \leq \trace(\bPhi_{V,n,t}^{\phi} \bX_t) \leq (V_n^{\mathrm{max}})^2, \forall\, \phi, \, \forall \, n \in \cN \label{mg2-voltage} \\
& \sum_{t \in \cI_{d,n}^{\phi}} \bar{P}_{d,n,t}^{\phi} \Delta_t = E_{d,n}^{\phi}, \hspace{.1cm}  \forall \,d \in \cD_n^{\phi}, \phi \in \cP_n, n \in \cN  \label{mg2-controllable} \\
&  \hspace{.6cm} 0 \leq \bar{P}_{d,n,t}^{\phi} \leq  \bar{P}_{d,n}^{\textrm{max}}, \hspace{.2cm} \forall \,\,  t \in \cI, d \in \cD_n^{\phi},\,  n \in \cN \label{mg2-Pdlimits} \\
& \rank(\bX_t) = 1,  \hspace{.2cm} \forall \,\,  t \in \cI   \label{mg2-rank} \\
& \bX_t \succeq \mathbf{0}. \,\, [\bX_t]_{\cP_0,\cP_0} = \mathbf{1}_{3\times3},  \hspace{.2cm} \forall \,\,  t \in \cI \, .
\end{align}
\end{subequations}
Unfortunately, (P2) is still nonconvex because of the rank-1 constraint on the positive semi-definite matrices $\{\bX_t\}$. Nevertheless, (P2) is amenable to the SDP relaxation technique, which amounts to dropping the rank constraints, thus relaxing nonconvex problems to SDP ones; see e.g., the tutorial~\cite{luospmag10}, and the works in~\cite{ZhuNAPS11} and~\cite{LavaeiLow}, where this technique is employed for power system state estimation and OPF for power transmission systems, respectively. Leveraging the SDP relaxation technique here too, it is possible to obtain the following \emph{convex} relaxation of (P2): 
\begin{subequations}
\begin{align} 
& \hspace{-.5cm} (P3) \,\,  \min_{\{\bX_t\},\cV_d} \sum_{t \in \cI} \kappa_{0,t} \sum_{\phi \in \cP_0} \trace(\bPhi_{P,0,t}^{\phi} \bX_t) \nonumber \\ 
& \hspace{1.1cm} + \sum_{t \in \cI}  \sum_{s \in \cS}  c_{s,t} \sum_{\phi \in \cP_s} \trace(\bPhi_{P,s,t}^{\phi} \bX_t)  \\
\textrm{s.t.} & \,\, \bX_t \succeq \mathbf{0},  \,\, [\bX_t]_{\cP_0,\cP_0} = \mathbf{1}_{3\times3}, \textrm{ and } \eqref{mg2-P}-\eqref{mg2-Pdlimits} \, . \nonumber
\end{align}
\end{subequations}
Clearly, if all the optimal matrices $\{\bX_{t}^{\textrm{opt}}\}$ of (P3) have rank $1$, then the variables $\{\bX_{t}^{\textrm{opt}}\}, \cV_{d}^{\textrm{opt}}$ represent a globally optimal solution also for (P2). Further, since (P1) and (P2) are \emph{equivalent}, there exist $2 |\cI|$ vectors $\{\bx_{t}^{\textrm{opt}}\}$ and $\{\bv_{t}^{\textrm{opt}}\}$, with $\bX_{t}^{\textrm{opt}} = \bx_{t}^{\textrm{opt}} \bx_{t}^{\textrm{opt} \cH}$ and $\bv_{t}^{\textrm{opt}} := \ba_{0,t} \circ \bx_{t}^{\textrm{opt}}$, for all $t\in \cI$, such that the optimal objective functions of (P1) and (P2) coincide. This is formally summarized next.  

\begin{proposition}  
Let $\{\bX_{t}^{\textrm{opt}}\}, \cV_{d}^{\textrm{opt}}$ be the optimal solution of (P3), and assume that $\rank(\bX_{t}^{\textrm{opt}}) = 1$, for all $t \in \cI$. Then, a globally optimal solution of (P1) is given by $\cV_{d}^{\textrm{opt}}$, the vectors of complex line-to-ground voltages
\begin{align} 
\bv_{t}^{\textrm{opt}} := \sqrt{\lambda_{1,t}} \bD_{0,t} \bu_{1,t} \, , \quad \forall \,\, t \in \cI
\end{align}
where $\lambda_{1,t} \in \mathbb{R}^+$ is the unique non-zero eigenvalue of $\bX_{t}^{\textrm{opt}}$ and $\bu_{1,t}$ the corresponding eigenvector, and the supplied active and reactive powers
\begin{align} 
P_{G,s,t}^{\textrm{opt}} & = \trace(\bPhi_{P,s,t}^{\phi} \bv_{t}^{\textrm{opt}} \bv_{t}^{\textrm{opt} \cH} ) + P_{L,s,t}^{\phi} + \sum_{d \in \cD_s^{\phi}} \bar{P}_{d,s,t}^{\phi, \textrm{opt}} \\
P_{Q,s,t}^{\textrm{opt}} & = \trace(\bPhi_{Q,s,t}^{\phi} \bv_{t}^{\textrm{opt}} \bv_{t}^{\textrm{opt} \cH} ) + Q_{L,s,t}^{\phi} \, , \,\forall s \in \cS \cup \{0\} .
\end{align}
\label{prop:rank1}
\end{proposition} 
%\hfill $\Box$

The upshot of the proposed formulation is that the \emph{globally} optimal solution of (P2) (and hence (P1)) can be obtained via standard interior-point solvers, in \emph{polynomial-time} complexity; see, for example, the complexity bounds for SDP reported in~\cite[Ch.~4]{Nemirovski_lecture} and~\cite{luospmag10}. This is in contrast with gradient descent-based solvers for nonconvex programs, sequential quadratic programming, and particle swarm optimization, which in general do not guarantee global optimality of the obtained solutions, face challenges pertaining to sensitivity of the initial guess, convergence, and complexity which increases with the number of iterations. Notice also, that matrices $\{\bPhi_{P,n,t}^{\phi},\bPhi_{Q,n,t}^{\phi}, \bPhi_{V,n}^{\phi}\}$ are very sparse. This property can be leveraged to substantially reduce the computational burden of interior-point solvers; for instance, the so-called ``chordal'' structure of matrices $\{\bX_t\}$ can be effectively exploited, as advocated in~\cite{Jabr12}.%\footnote{Leveraging these special structures to develop interior-point iterations tailored for (P3) goes beyond the scope and space limits of this paper}%, but constitutes an interesting future research direction.} 

Since (P3) is a relaxed version of (P2), matrices $\bX_{t}^{\textrm{opt}}$ could have rank greater than $1$. In this case, rank reduction techniques can be employed to find a feasible rank-1 approximation of $\bX_{t}^{\textrm{opt}}$ (see~\cite{luospmag10} and references therein). The resultant solution is feasible for (P2), but generally suboptimal~\cite{luospmag10}. Notably, when \emph{balanced} distribution networks are considered,~\cite{Tse12} established conditions on the voltage angles and the reactive power injections under which rank-$1$ matrices are \emph{always} obtained provided the non-relaxed problem is feasible. 
%The proof amounts to showing that the Pareto-front remains unchanged upon taking the convex hull of the overall power injection region, which is given by the product of the ellipsoidal line flow regions. 
Derivation of similar conditions in the present context constitutes an interesting future research direction, that will naturally complement the result in~\cite{Tse12}.

%%%%%%%%%%%%%%%%%%%%%%%%%%%%%
\subsection{Feeders with two- and single-phase lines}
\label{sec:mixed}
%%%%%%%%%%%%%%%%%%%%%%%%%%%%%

For feeders with two- and single-phase laterals and sub-laterals, the dimensions of matrix $\bY$ have to be adjusted to $\sum_{n=0}^N |\cP_n| \times \sum_{n=0}^N |\cP_n|$, and its entries have to be as follows: 

\emph{i)} matrix $-\bZ_{nm}^{-1}$ occupies the $|\cP_{mn}| \times |\cP_{mn}|$ off-diagonal block corresponding to line $(m,n) \in \cE$; and, 

\emph{ii)} the $|\cP_{n}| \times |\cP_{n}|$ diagonal block corresponding to node $n \in \cN \cup \{0\}$ is obtained as   
\begin{align} 
[\bY]_{\cP_n,\cP_n} := \sum_{m \in \cN_m}  \left(\frac{1}{2}\tilde{\bY}_{mn}^{(s)} + \tilde{\bZ}_{mn}^{-1} \right)
\end{align}
where $\tilde{\bZ}_{mn} = \bZ_{mn}$ and $\tilde{\bY}_{mn}^{(s)} = \bY_{mn}^{(s)}$ if $\cP_n = \cP_{mn}$, otherwise $[\tilde{\bZ}_{mn}]_{\cP_{nm},\cP_{nm}} = \bZ_{mn}$ and $[\tilde{\bZ}_{mn}]_{\cP_n \backslash \cP_{nm},\cP_n \backslash \cP_{nm}} = 0$ ($\tilde{\bY}_{mn}^{(s)}$ is computed likewise). Re-defining the $\sum_{n=0}^N |\cP_n| \times 1$ vectors collecting voltages and currents as $\bv^\prime_t := [\bv_{0,t}^\cT, [\bv_{1,t}]_{\cP_{1}}^\cT,\ldots, [\bv_{N,t}]_{\cP_{N}}]^T$ and $\bi^\prime_t := [\bi_{0,t}^\cT, [\bi_{1,t}]_{\cP_{1}}^\cT,\ldots, [\bi_{N,t}]_{\cP_{N}}^\cT]^\cT$, respectively,~\eqref{currents} can be re-written again in vector-matrix form as $\bi^\prime_t = \bY \bv^\prime_t$, for all $t \in \cI$, and the procedure~\eqref{eq:Phi}--\eqref{mg2} can be readily followed. 

%Clearly, constraints corresponding to the missing phases on the buses have to be discarded in (P2) and (P3).

  %%%%%%%%%%%%%%%%%%%%%%%%%%%%%%%%%%%%%%%%%%%%%%%
\section{Feasible voltage profile}
\label{sec:feasibility}
%%%%%%%%%%%%%%%%%%%%%%%%%%%%%%%%%%%%%%%%%%%%%%%

%Since the high resistance-reactance ratio in conventional distribution lines renders the voltages rather sensitive to variations in the net powers injected at the nodes, sections of feeders may face issues pertaining to over- and under-voltage conditions when elastic loads and DG units are present. Therefore, 

To effect voltage regulation, and avoid  abrupt voltage drops, constraints on the minimum and maximum utilization and service voltages were imposed in (P1).  Constraints~\eqref{mg-Vlimits} however, may challenge the feasibility of (P1), since it may not be possible to meet the minimum (maximum) utilization and service voltage requirements when feeders are heavily stressed and DG units supply a substantial amount of power. It is thus of prime importance to perform preemptive analysis of the feasible voltage profile in order to unveil possible infeasibility of (P1) and, in case, facilitate corrective actions. To this end, consider solving the following optimization problem
\begin{subequations}
\label{Pmg4}
\begin{align} 
& \hspace{-.5cm} (P4)\,\, \min_{\substack{\cV^{(1)}, \cV_d }} \,\, (1-w_V) \sum_{t \in \cI}  \sum_{n \in \cN}  \sum_{\phi \in \cP_n} \left( |V_{n,t}^{\phi}|^2 - |V_{n}^{\textrm{ref}}|^2 \right)^2 \nonumber \\
& + w_V \sum_{t \in \cI}  \left(\kappa_{0,t} \sum_{\phi \in \cP_0} P_{0,t}^{\phi} + \sum_{s \in \cS}  c_{s,t} \sum_{\phi \in \cP_s} P_{G,s,t}^{\phi} \right) \label{mg-cost4} \\
& \textrm{s.t.} \,\, \eqref{mg-balance-source}-\eqref{mg-Pdlimits}, \textrm{ and }~\eqref{mg-plimits}-\eqref{mg-qlimits} \nonumber
\end{align}
\end{subequations}
with $w_V \in (0,1)$ denoting a weighting coefficient, and $|V_{n}^{\textrm{ref}}|$ the prescribed voltage magnitude of the feeder (e.g., $|V_{n}^{\textrm{ref}}| = 1$ p.u.). Although constraints~\eqref{mg-Vlimits} are not present in (P4), the first term in~\eqref{mg-cost4}  promotes regulation by penalizing voltage magnitudes that deviate from the nominal ones.

Similar to (P1), problem~\eqref{Pmg4} is nonconvex. However, by exploiting again Lemma~\ref{SDPreformulation}, along with the SDP relaxation technique, the following convex problem is obtained
\begin{subequations}
\label{Pmg5}
\begin{align} 
& \hspace{-.5cm} (P5)\,\, \min_{\substack{\{\bX_t\}, \cV_d \\ \{\alpha_{n,t}^{\phi} \}}} \,\, (1-w_V) \sum_{t, n, \phi}  \alpha_{n,t}^{\phi}  \nonumber \\ 
& \hspace{1.5cm} + w_V \sum_{t \in \cI} \kappa_{0,t} \sum_{\phi \in \cP_0} \trace(\bPhi_{P,0,t}^{\phi} \bX_t) \nonumber \\ 
& \hspace{1.5cm} + w_V \sum_{t \in \cI}  \sum_{s \in \cS}  c_{s,t} \sum_{\phi \in \cP_s} \trace(\bPhi_{P,s,t}^{\phi} \bX_t) \label{mg-cost5} \\
\textrm{s.t.} &  \nonumber \\
& \left[ \begin{array}{ll}
- \alpha_{n,t}^{\phi} & \trace(\bPhi_{V,n}^{\phi} \bX_t)-|V_{n}^{\textrm{ref}}|^2 \\
\trace(\bPhi_{V,n}^{\phi} \bX_t)-|V_{n}^{\textrm{ref}}|^2 & -1
\end{array}
\right] \preceq \mathbf{0} \label{mg5-alpha} \\
& \,\, \bX_t \succeq \mathbf{0},  \,\, [\bX_t]_{\cP_0,\cP_0} = \mathbf{1}_{3\times3}, \textrm{ and } \eqref{mg2-P}-\eqref{mg2-Pdlimits}  \nonumber
\end{align}
\end{subequations}
where constraint~\eqref{mg5-alpha} is enforced for all nodes $n$, per phase $\phi$, and time slot $t$. Notice that by using~\eqref{eq:V} the first term in~\eqref{mg-cost4} becomes quadratic in $\{\bX_t\}$. To bypass this hurdle, the non-negative real variables $\{\alpha_{n,t}^{\phi} \}$ are introduced to upper bound each term $(\trace(\bPhi_{V,n}^{\phi} \bX_t)-|V_{n}^{\textrm{ref}}|^2)^2$, and   the Schur's complement is subsequently employed to obtain~\eqref{mg5-alpha}. If for a given $w_V$, all matrices $\{\bX_t\}$ have rank $1$, then the optimal solution of (P5) is also a globally optimal solution of (P4) (see Proposition~\ref{prop:rank1}).  

Clearly, if the voltages $\{V_{n,t}^{\phi, \mathrm{opt}}\}$ obtained from (P5) satisfy $(V_{n,t}^{\mathrm{min}})^2 \leq |V_{n,t}^{\phi, \mathrm{opt}}|^2 \leq (V_{n,t}^{\mathrm{max}})^2$, for all $t, n, \phi$, then it is possible to proceed with the solution of the economic dispatch problem (P3). On the other hand, if some of the voltage magnitudes largely deviate from $|V_{n}^{\textrm{ref}}|$, corrective actions have to be taken; these include, for example, switching the taps of controllable capacitor banks, or curtailing portion(s) of the loads.

%%%%%%%%%%%%%%%%%%%%%%%%%%%%%%%%%%%%%%%%%%%%%%%
\section{Thermal and Quality-of-power constraints}
\label{sec:Additional}
%%%%%%%%%%%%%%%%%%%%%%%%%%%%%%%%%%%%%%%%%%%%%%%

%%%%%%%%%%%%%%%%%%%%%%%
\subsection{Thermal constraints}
\label{sec:linecurrents}
%%%%%%%%%%%%%%%%%%%%%%%

High current magnitudes on the lines can have detrimental effects on both efficiency and reliability of the distribution network. From an economical perspective, an additional (real) power has to be drawn from the main grid, or, supplied by the DG units in order to compensate for the increased power dissipated on the distribution lines. On the other hand, conductors may overheat if stressed by high currents over a prolonged time interval, and may eventually fail. This, in turn, would trigger an outage event, with consequent interruption of the power delivery in portions of the feeder. 

To alleviate these concerns, it is of interest to constrain either the power dispelled on the conductors, or, the magnitude of currents flowing through the distribution lines~\cite{Tse12}, which amounts to adding one of the following constraints in (P1):
\begin{align}
|I_{mn,t}^{\phi}| & \leq  I_{mn}^{\textrm{max}}     \label{current} \\
P_{mn}^{\phi}  :=  |I_{mn}^{\phi}|^2 \Re\{[\bZ_{mn}]_{\{\phi\},\{\phi\}}\}   & \leq  P_{mn}^{\textrm{max}}       \label{powerdissipated}
\end{align}
where $I_{mn,t}^{\phi}$ and $P_{mn}^{\phi}$ denote the current flowing on the phase $\phi$ of line $(m,n) \in \cE$, and the active power lost on the same line and phase, respectively. 

An SDP-consistent re-formulation of~\eqref{current}--\eqref{powerdissipated} has to be derived in order to accommodate the aforementioned constraints in (P2) and (P3). To this end, let $\bi_{mn,t} : = [\{I_{mn,t}^{\phi}\}]^\cT$ denote the $|\cP_{mn}| \times 1$ vector collecting the complex currents flowing through line $(m,n) \in \cE$, which is related to the line-to-ground voltages $\bv_{n,t}$ and $\bv_{n,t}$ as (cf.~\eqref{currents}) 
\begin{equation}
\bi_{mn,t}  =  \bZ_{mn}^{-1} \left([\bv_{m,t}]_{\cP_{mn}} - [\bv_{n,t}]_{\cP_{mn}} \right)  \, .    \label{linecurrent}
\end{equation}
Notice that since $\bZ_{mn}^{-1}$ is generally not diagonal~\cite{testfeeder},~\eqref{linecurrent} captures also current components arising from mutual inductive reactances and capacitive coefficients. Define the $|\cP_{mn}| \times \sum_{n=0}^N |\cP_n|$ complex matrix
\begin{align}
\bB_{mn} & := [\mathbf{0}_{|\cP_{mn}| \times \sum_{n=0}^{m-1} |\cP_n|}, \check{\bZ}_{mn}^m, \ldots \nonumber \\
& \mathbf{0}_{|\cP_{mn}| \times \sum_{n=m+1}^{n-1} |\cP_n|}, \check{\bZ}_{mn}^n  \mathbf{0}_{|\cP_{mn}| \times \sum_{n=n+1}^{N} |\cP_n|}]
\end{align}
where $\check{\bZ}_{mn}^m$ is a $|\cP_{mn}| \times |\cP_{m}|$ matrix with elements $[\check{\bZ}_{mn}^m]_{\cP_{mn},\cP_{mn}} = \bZ_{mn}^{-1}$ and $[\check{\bZ}_{mn}^m]_{\cP_{mn},\cP_m \backslash \cP_{mn}} = \mathbf{0}$; likewise, $\check{\bZ}_{mn}^n$ has dimensions $|\cP_{mn}| \times |\cP_{n}|$, and its entries are filled as $[\check{\bZ}_{mn}^n]_{\cP_{mn},\cP_{mn}} = -\bZ_{mn}^{-1}$ and $[\check{\bZ}_{mn}^n]_{\cP_{mn},\cP_n \backslash \cP_{mn}} = \mathbf{0}$. Thus, building on~\eqref{linecurrent}, an SDP-compliant re-formulation of~\eqref{current}--\eqref{powerdissipated} is provided in the following lemma.     

\vspace{2mm}

\begin{lemma}
\label{limitcurrent}
Consider the Hermitian matrix 
\begin{align}
\bPhi_{I,mn,t}^{\phi} :=  \bD_{0,t}^{\cH} \bB_{mn}^{\cH} \be_{mn}^{\phi} (\be_{mn}^{\phi})^{\cT} \bB_{mn} \bD_{0,t}
\end{align}
where $\{\be_{mn}^\phi\}_{\phi \in \cP_{mn}}$ denotes the canonical basis of $\mathbb{R}^{|\cP_{mn}|}$. Then, constraints~\eqref{current}--\eqref{powerdissipated} can be expressed linearly in the outer-product $\bX_t$ as   
\begin{align}
\trace\{ \bPhi_{I,mn,t}^{\phi} \bX_t \} & \leq  I_{mn}^{\textrm{max}}  \label{currentSDP} \\
\Re\{\trace\{ \be_{mn}^{\phi} (\be_{mn}^{\phi})^{\cT}  \bZ_{mn}\} \}  \trace\{ \bPhi_{I,mn,t}^{\phi} \bX_t \} & \leq  P_{mn}^{\textrm{max}}.   \label{currentSDP} 
\end{align}
\end{lemma}
\emph{Proof.}   See the Appendix.  \hfill $\Box$

\vspace{2mm}

In distribution feeders, line outage events maybe triggered by overheating effects on the neutral cable(s), especially those experiencing highly unbalanced load conditions. Towards deriving constraints on the magnitude of neutral current(s), let $\cP_{mn}^{(\varphi)}$ denote the set of grounded neutral cables that are present on the line $(m,n) \in \cE$, and $\bT_{mn}$ the  $|\cP_{mn}^{(\varphi)}| \times |\cP_{mn}|$ neutral transformation matrix obtained from the primitive impedance matrix of the distribution line via Kron reduction~\cite[Sec.~4.1]{Kerstingbook}. For example, the neutral transformation matrix of a four-wire grounded wye segment has dimensions $1 \times 3$, while its dimensions increase to $3 \times 3$ for an underground wye line with three neutral conductors. Thus, the neutral currents $\bi^{(\varphi)}_{mn,t} := [I^{(1)}_{mn,t}, \ldots, I^{(N^{\varphi})}_{mn,t}]^{\cT}$ are linearly related to the line currents $\bi_{mn,t}$ as~\cite[Sec.~4.1]{Kerstingbook} 
\begin{equation}
\bi^{(\varphi)}_{mn,t} = \bT_{mn} \bi_{mn,t} \, .     \label{neutralcurrent}
\end{equation}
It readily follows from~\eqref{neutralcurrent} and the result of Lemma~\ref{limitcurrent}, that the magnitude of the current on the neutral cables can be constrained in the SDP problem (P3) as  
\begin{equation}
\trace\{ \bPhi_{I,mn,t}^{(\varphi)} \bX_t \}  \leq  I_{mn}^{(\varphi), \textrm{max}},    \,\, \forall \, \varphi \in \cP^{(\varphi)}_{mn}   \label{neutralcurrentlimit}
\end{equation}
with 
\begin{align}
\bPhi_{I,mn,t}^{(\varphi)} :=  \bD_{0,t}^{\cH} \bB_{mn}^{\cH} \bT_{mn}^{\cH} \be_{mn}^{(\varphi)} (\be_{mn}^{(\varphi)})^{\cT} \bT_{mn} \bB_{mn} \bD_{0,t}
\end{align}
where, as usual, $\{\be_{mn}^{(\varphi)}\}$ is the canonical basis of $\mathbb{R}^{|\cP^{(\varphi)}_{mn}|}$.

There is an increasing concern on the thermal effects arising from harmonic currents in the neutral cable(s). In this case, constraints similar to~\eqref{neutralcurrentlimit} can be imposed on a per-harmonic basis (see, e.g.~\cite{Forner}).

%%%%%%%%%%%%%%%%%%%%%%%
\subsection{Constraints on the power factor}
\label{sec:pfcorrection}
%%%%%%%%%%%%%%%%%%%%%%%

The PF has been increasingly recognized as one of the principal measures of efficiency and reliability of power distribution networks~\cite{Momoh09,Roytelman00,Forner}. High PF translates to lower generation and transmission costs, and enhanced protection of transmission lines from overheating (hence, higher resilience to line outages). Constraining the PF at the PCC is tantamount to limiting the reactive power exchanged with the main power grid. This, in turn, has two well-appreciated merits: \emph{i)} it alleviates the power losses experienced along the backbone of the feeder~\cite{Roytelman00}; and, \emph{ii)} it limits the current drawn at the PCC, and therefore facilitates coexistence of multiple feeders on the same distribution line or substation without requiring components such as, e.g. conductors, transformers, and switchgear of increased size. 

Unfortunately, the definition of PF for an unbalanced polyphase system is not unique~\cite{Emanuel93}. In this paper, a per-phase definition is adopted in order to limit the reactive power exchanged at the PCC on each phase. Intuitively, polyphase variants~\cite{Emanuel93} may induce high discrepancies between the amount of reactive power exchanged per phase, but with a ``good'' polyphase PF nevertheless. Let $\eta_{0,t}^{\phi} \in [0,1]$ denote the minimum PF required at the PCC on the phase $\phi \in \cP_0$ at time slot $t \in \cI$. Now consider adding the following constraint to (P1): 
\begin{equation}
P_{0,t}^{\phi} \left(|V_{0,t}^{\phi}| |I_{0,t}^{\phi}| \right)^{-1}  \geq \eta_{0,t}^{\phi}, \quad \quad \forall \,\, t \in \cI, \phi \in \cP_0 \label{mg-powerfactor}
\end{equation}
where voltages $\{|V_{0,t}^{\phi}| \}$ are known~\cite{Kerstingbook}. Notice that DG units complement the power supplied by the utility, and are usually not sufficient to satisfy the load demand on their own. Under this premise, an SDP-consistent reformulation of~\eqref{mg-powerfactor} can be readily obtained, as summarized next.

\vspace{.2cm}

\begin{lemma}  
Provided the power supplied by the DG units does not exceed the total load demand at the feeder,~\eqref{mg-powerfactor} is equivalently expressed as a linear function of $\bX$ as 
\begin{align}
\left \{\begin{array}{ll} 
\tilde{\eta}_{0,t}^{\phi}  \trace(\bPhi_{P,n}^{\phi} \bX_t)  -  \trace(\bPhi_{Q,n}^{\phi} \bX_t)  \geq 0 \\
\tilde{\eta}_{0,t}^{\phi}  \trace(\bPhi_{P,n}^{\phi} \bX_t)  +  \trace(\bPhi_{Q,n}^{\phi} \bX_t)  \geq 0 \, .
\end{array} \right.
\end{align}
\end{lemma}
\hfill $\Box$ 

\vspace{.2cm}

Additional charges are generally applied to residential and industrial loads with a poor PF. In the presence of highly inductive loads, 
capacitor banks are usually employed to balance reactive demand, and thus maintain the PF as close as possible to 1~\cite{Kerstingbook}. Recall that $y_{C,n}^{\phi}$ denotes the susceptance of a capacitor connected at node $n$ and phase $\phi$, and the provided reactive power amounts to $Q_{C,n}^{\phi} = y_{C,n}^{\phi} |V_{n,t}^{\phi}|^2$. With $Q_{L,n,t}^{\phi} > 0$ denoting the reactive power demanded by an inductive load, a minimum per-phase PF $\eta_{n,t}^{\phi} \in [0,1]$ is imposed as [cf.~\eqref{mg-powerfactor}] 
\begin{align} 
\label{powerfactor_node}
\frac{P_{L,n,t}^{\phi}}{\sqrt{(P_{L,n,t}^{\phi})^2 + (Q_{L,n,t}^{\phi} - Q_{C,n,t}^{\phi})^2 }} \geq \eta_{n,t}^{\phi} \, \quad \forall \, \phi \in \cP_n
\end{align}
where $P_{L,n,t}^{\phi}$ is given. Clearly, $|V_{n,t}^{\phi}|^2$ can be re-expressed as a linear function of $\bX_t$ using~\eqref{eq:V}, and~\eqref{powerfactor_node} can be reformulated to obtain the following SDP-compliant form. 

\vspace{.2cm}

\begin{lemma}
\label{SDPformPF}
Using~\eqref{eq:V}, constraint~\eqref{powerfactor_node} is equivalent to the following linear matrix inequality (with $Q_{C,n,t}^{\phi} = y_{C,n}^{\phi}\trace(\bPhi_{V,n}^{\phi} \bX_t)$)
\begin{align} 
\label{eq:SDPformPF}
\left[ \begin{array}{ll}
- \left(\frac{P_{L,n,t}^{\phi}}{\eta_{n,t}^{\phi}} \right)^2 & P_{L,n,t}^{\phi} + Q_{L,n,t}^{\phi} - Q_{C,n,t}^{\phi}   \\ 
P_{L,n,t}^{\phi} + Q_{L,n,t}^{\phi} - Q_{C,n,t}^{\phi} & -1
\end{array}
\right] \preceq \mathbf{0}.
\end{align} 
\end{lemma}
\emph{Proof.}   See the Appendix.  \hfill $\Box$

Controllable capacitor banks can be accounted for by associating an integer variable with each of the capacitor switches. To tackle the resultant mixed integer nonlinear problem, exhaustive search over the switches can be performed~\cite{Paudyal11}. In this case, (P3) is solved for each switch configuration.

%%%%%%%%%%%%%%%%%%%%%%%%%%%%%%%%%%%%%%%%%%%%%%%
\section{Numerical Tests}
\label{sec:results}
%%%%%%%%%%%%%%%%%%%%%%%%%%%%%%%%%%%%%%%%%%%%%%%

The proposed optimization framework for unbalanced three-phase systems is tested on the IEEE 13-node test feeder shown in Fig.~\ref{fig:feeder}. Compared to the original scheme in~\cite{testfeeder}, DG units are placed at nodes 1 and 10. Specifically, single-phase DG units supply a maximum real power of $300$ kW and $500$ kW, respectively, and they operate at a unitary PF. Capacitor banks with rated reactive power of $200$ kVAr and $100$ kVAr are present at nodes $5$ and $8$, respectively. Line impedance and shunt admittance matrices are computed based on the dataset in~\cite{testfeeder}.
To solve (P3) (and (P5) for a preemptive voltage profile analysis), the MATLAB-based optimization modeling package \texttt{CVX}~\cite{cvx} is used, along with the interior-point based solver \texttt{SeDuMi}~\cite{SeDuMi}. 

%%%%%%%%%%%%%%%%%%%%
\begin{figure}[t]
\begin{center}
\includegraphics[width=0.35\textwidth]{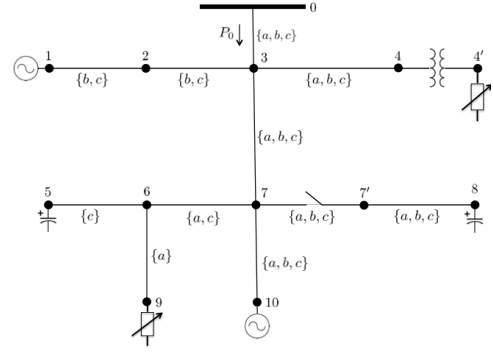}
\caption{Modified IEEE 13-bus test feeder.}
\label{fig:feeder}
\end{center}
\vspace{-.5cm}
\end{figure}
%%%%%%%%%%%%%%%%%%%%

The time horizon is $24$h, and slots of $1$h are considered; that is, $\cI = \{1\, \textrm{AM}, 2\, \textrm{AM}, \ldots, 11\, \textrm{PM}, 12\, \textrm{AM}\}$. The loads specified in~\cite{testfeeder} are assumed to be the peak demands of the day, and the ``spring mid-week'' load profiles reported in~\cite{report_Canada} are used to generate $\{P_{L,n,t}^{\phi}, Q_{L,n,t}^{\phi} \}_{t\in \cI}$. Specifically, the ``commercial load profile'' in~\cite[Sec.~1.1]{report_Canada} is used for node $9$, whereas the ``residential load profile'' is applied to all the remaining nodes; a Gaussian random variable with mean 1 and standard deviation $0.1$ is used to insinuate a perturbation on the profiles on a per-node and a per-time basis. The resulting aggregate real loads per phase are depicted in Fig.~\ref{fig:Profile} (similar trends are obtained for the reactive loads, but are not reported due to space limitations).  

Ten controllable loads are present at node $5^{\prime}$, and are allocated as follows: $2$ on phase $a$, $4$ on phase $b$, and $4$ on phase $c$. The energy requirement is $11$ kWh, and the cap $\bar{P}_{d,5,t}^{\textrm{max}}$ is set to $4$ kW, so as to resemble the demands of $10$ PHEVs~\cite{Driesen10}. Customers are assumed to plug-in the PHEVs at 6 PM, and the charging has to be completed by 6 AM. Two additional elastic loads are present at node $9$, and have to be satisfied between 8 AM and 4 PM.  In this case, the energy requirement per each load is $30$ kWh, and no cap is present for $\{\bar{P}_{d,9,t}^{\phi}\}$.

To model the price of the power purchased from the main distribution grid $\{\kappa_{0,t}\}$, one-day ahead locational marginal prices (LMPs) available in the Midwest Independent Transmission System Operator (MISO)~\cite{midwestiso} database are utilized. Specifically, the LMPs for the Minneapolis area on June 7th, 2012 are utilized throughout this section, and are reported in Fig.~\ref{fig:Profile}.  On the other hand, the cost incurred by the use of the DG units is kept constant over time, and it is set to $30$ $\$$/MW. 

A minimum PF of $\eta_{0,t}^{\phi} = 0.8$ is required at the PCC, and the limits $V_n^{\mathrm{min}} = 0.95$ p.u. and $V_n^{\mathrm{max}} = 1.05$ p.u. are imposed to enforce voltage regulation. Finally, a balanced flat voltage profile is assumed at the PCC, with $|V_{0,t}^{\phi}| = 1.02$ p.u., and $\angle V_{0,t}^{a} = 0^{\circ}$, $\angle V_{0,t}^{b} = 120^{\circ}$, and $\angle V_{0,t}^{c} = -120^{\circ}$. Constraints on the line currents are not considered, since this datum is not available in~\cite{testfeeder}. 

%%%%%%%%%%%%%%%%%%%%
\begin{figure}[t]
\begin{center}
\includegraphics[width=0.5\textwidth]{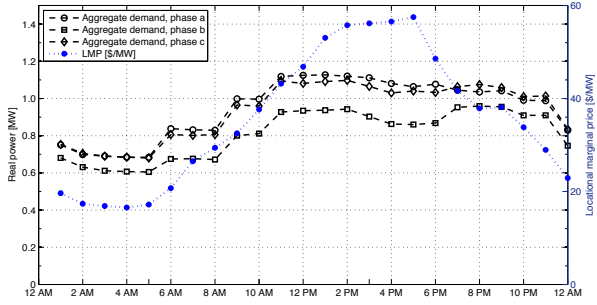}
\caption{Aggregate real load profile [MW], and prices $\{\kappa_{0,t}\}$ [$\$$/MW].}
\label{fig:Profile}
\end{center}
\vspace{-.6cm}
\end{figure}
%%%%%%%%%%%%%%%%%%%%

%%%%%%%%%%%%%%%%%%%%
\begin{figure}[t]
\begin{center}
\includegraphics[width=0.5\textwidth]{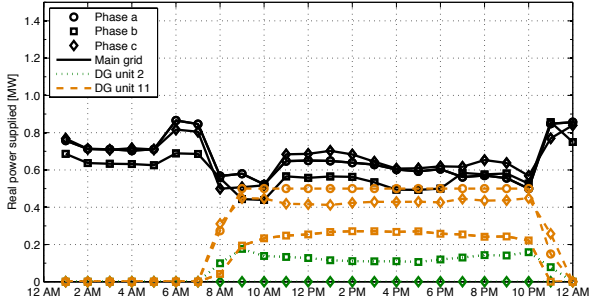}
\caption{Real power supplied [MW].}
\label{fig:supplied}
\end{center}
\vspace{-.5cm}
\end{figure}
%%%%%%%%%%%%%%%%%%%%

%%%%%%%%%%%%%%%%%%%%
\begin{figure}[t]
\begin{center}
\includegraphics[width=0.5\textwidth]{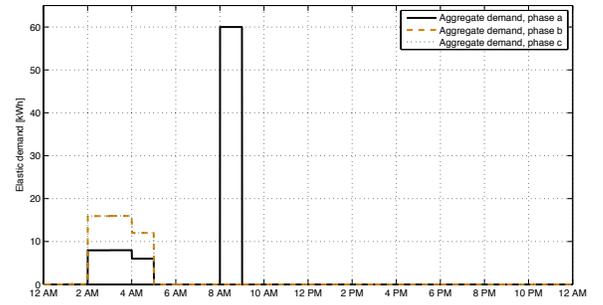}
\caption{Aggregate per-phase elastic demand [kWh].}
\label{fig:Elastic}
\end{center}
\vspace{-.5cm}
\end{figure}
%%%%%%%%%%%%%%%%%%%%

%%%%%%%%%%%%%%%%%%%%
\begin{figure}[t]
\begin{center}
\includegraphics[width=0.5\textwidth]{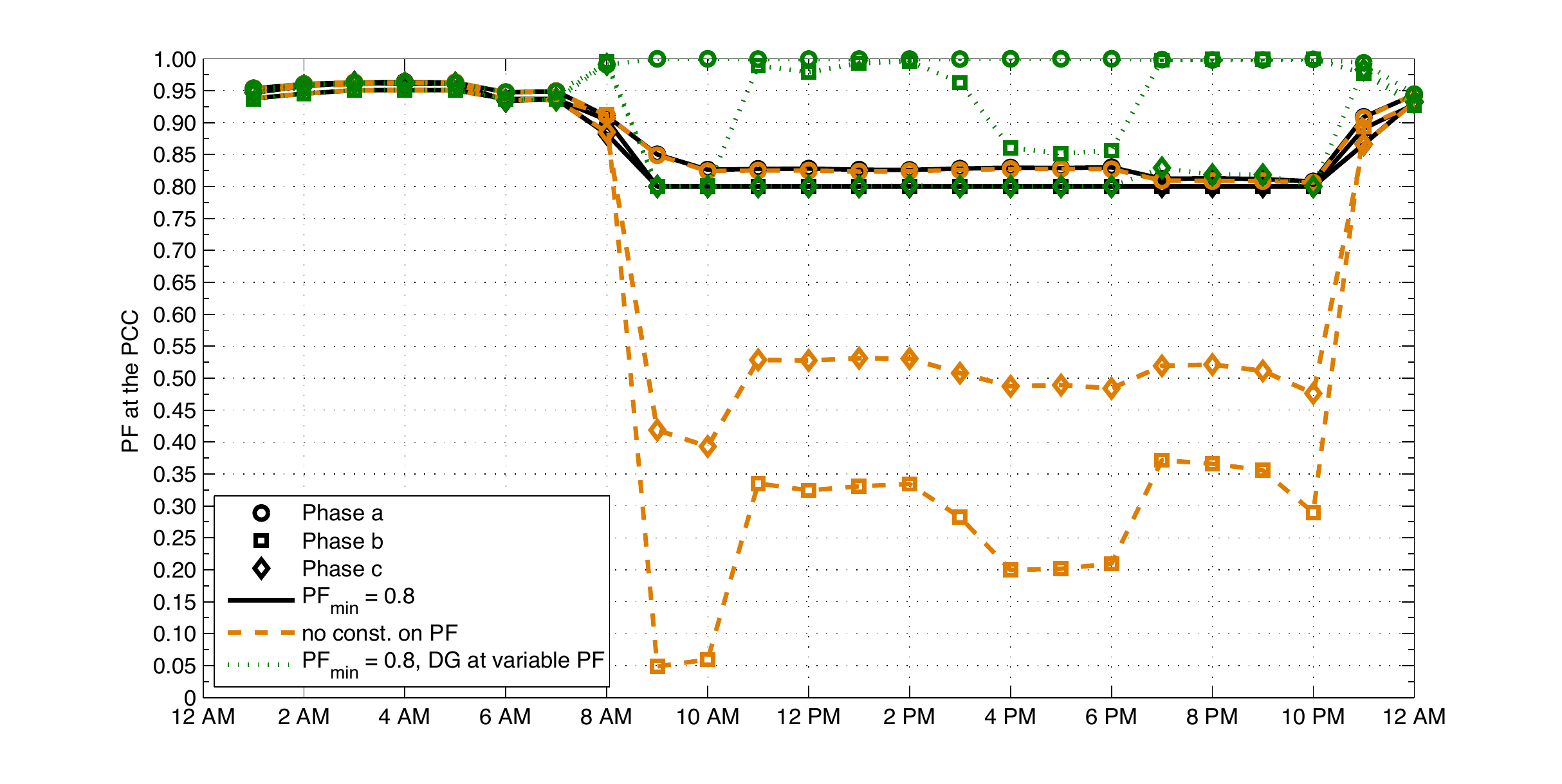}
\caption{PF at the substation.}
\label{fig:PFprofile}
\end{center}
\vspace{-.5cm}
\end{figure}
%%%%%%%%%%%%%%%%%%%%

Before proceeding, it is worth mentioning that the rank of the matrices $\{\bX_t^{\textrm{opt}}\}$ was \emph{always} 1 in the experiments reported in this section. Therefore, the \emph{globally} optimal solutions of  (P1) were always attained. This illustrates clearly the merits of the proposed formulation. 

Fig.~\ref{fig:supplied} depicts the active power supplied by the DG units, and drawn from the main distribution grid. As expected, the DG units are heavily utilized from $8$ AM to $11$ PM, which is the interval where the price of power purchased from the main distribution grid is higher that $30$ $\$$/MW. This, in turn, has the benefit of reducing the overall demand of the feeder during the peak hours (peak shaving). Notice however, that the DG units are not utilized at the maximum extent because of the constraint on the PF at the PCC~\cite{Kroposki10}, as it will be shown later on. The optimal allocation of the elastic energy demands is shown in Fig.~\ref{fig:Elastic}. It is shown that the PHEVs connected at node $5^\prime$ are charged from $2$ AM to $5$ AM, interval where both the LMPs and the non-deferrable loads are the lowest. A similar behavior is noticed for the elastic demands at node $9$; in fact, they are entirely satisfied in the time slot $[8 \, \textrm{AM}, 9 \, \textrm{AM}]$, which is the slot with the lowest LMP in the interval $[8 \, \textrm{AM}, 4 \, \textrm{PM}]$. 

Finally, Fig.~\ref{fig:PFprofile} portrays the trajectories of the PF at the PCC. It can be seen that the lower bound on the PF is tightly met when the DG units are active. In fact, as they supply real power to a lagging power system, a reduction of the PF is inevitably experienced at the PCC. This can be further noticed from the dotted (orange) trajectories, which correspond to the case where (P3) is solved without the constraints on the PF. In this case, the majority of the real power is supplied by the DG units, thus entailing a significant drop of the PF at the PCC. The case where the DG units can operate at a variable PF, from $0.5$ to $1$, is also considered. In this case, the DG units can supply a sufficiently amount of reactive power, so that the PF at the substation can be kept close to the unity most of the times. This suggests that DG units can be effectively utilized for providing reactive support~\cite{Kroposki10}, although an appropriate modeling of the cost incurred in this case is required.

%%%%%%%%%%%%%%%%%%%%%%%%%%%%%%%%%%%%%%%%%%%%%%%
\section{Concluding remarks}
\label{sec:conclusions}
%%%%%%%%%%%%%%%%%%%%%%%%%%%%%%%%%%%%%%%%%%%%%%%

The paper considered the economic dispatch problem for unbalanced three-phase power distribution networks, where the costs of power drawn from the main grid and supplied by the DG units over a given time horizon was minimized, while meeting the overall load demand and effecting voltage regulation. Is spite of the inherent non-convexity of the formulated problem, the SDP relaxation technique was advocated to obtain a (relaxed) convex problem. As corroborated by numerical tests, the main merit of the 
proposed approach consists in offering the potential of finding the globally optimal solution of the original nonconvex economic dispatch problem.

%Future efforts include consideration of controllable capacitor banks, voltage regulators, and renewable-based DG units. 

%%%%%%%%%%%%%%%%%%%%%%%%%%%%%%%%%%%%%%%%%%%%%%
\appendix
%%%%%%%%%%%%%%%%%%%%%%%%%%%%%%%%%%%%%%%%%%%%%%

\emph{Proof of Lemma}~\ref{SDPreformulation}. To prove~\eqref{eq:P}, notice first that the injected apparent power at node $n$, phase $\phi$ and time $t$ is given by $V_{n,t}^{\phi} (I_{n,t}^{\phi})^* = (V_{n,t}^{\phi, *} I_{n,t}^{\phi})^* = (\bv_t^\cH  \bar{\be}_n^{\phi} (\bar{\be}_n^{\phi})^{\cT} \bi_t)^{\cH}$. Next, noticing that $ \bv_t = \ba_{0,t} \circ \bx_t = \bD_{0,t} \bx_t$ and using $\bi_t = \bY \bv_t$, it follows that $(\bv_t^\cH  \bar{\be}_n^{\phi} (\bar{\be}_n^{\phi})^{\cT} \bi_t)^{\cH} = (\bx_t^\cH \bD_{0,t}^\cH  \bar{\be}_n^{\phi} (\bar{\be}_n^{\phi})^{\cT} \bY \bD_{0,t} \bx_t)^{\cH} = (\bx_t^\cH \bD_{0,t}^\cH  \bY^{\phi}_n \bD_{0,t} \bx_t)^{\cH} = \bx_t^\cH \bD_{0,t}^\cH  (\bY^{\phi}_n)^{\cH} \bD_{0,t} \bx_t$, which can be equivalently rewritten as $\trace(\bD_{0,t}^\cH \bY_n^{\phi} \bD_{0,t} \bX_t)$. Thus, the injected real and reactive powers can be obtained by using, respectively, the real and imaginary parts of $(\bY^{\phi}_n)^{\cH}$. Finally,~\eqref{eq:V} can be readily established by noticing that $|V_{n,t}^{\phi}|^2 = \bv_t^\cH \bar{\be}_n^{\phi} (\bar{\be}_n^{\phi})^{\cT} \bv_t = \bx_t^\cH \cD_{0,t}^\cH \bar{\be}_n^{\phi} (\bar{\be}_n^{\phi})^{\cT} \bD_{0,t} \bx_t = \trace(\bD_{0,t}^\cH \bar{\be}_n^{\phi} (\bar{\be}_n^{\phi})^{\cT} \bD_{0,t} \bX_t)$. 

\emph{Proof of Lemma}~\ref{limitcurrent}. From~\eqref{linecurrent}, and using the definitions of $\check{\bZ}_{mn}^m$ and $\check{\bZ}_{mn}^n$, it follows that $\bi_{mn,t} \bi_{mn,t}^{\cH} = \check{\bZ}_{mn}^m \bv_{m,t} \bv_{m,t}^{\cH}  \check{\bZ}_{mn}^{m \cH} + \check{\bZ}_{mn}^n \bv_{n,t} \bv_{n,t}^{\cH}  \check{\bZ}_{mn}^{n \cH} - \check{\bZ}_{mn}^m \bv_{m,t} \bv_{n,t}^{\cH}  \check{\bZ}_{mn}^{n \cH} - \check{\bZ}_{mn}^n \bv_{n,t} \bv_{m,t}^{\cH}  \check{\bZ}_{mn}^{m \cH} = \bB_{mn} \bv_{t} \bv_{t}^{\cH}   \bB_{mn}^{\cH}$. Thus, $|I_{mn}^{\phi}|^2$ is given by $|I_{mn}^{\phi}|^2 = (\be_{mn}^{\phi})^{\cT} \bB_{mn} \bD_{0,t} \bx_{t} \bx_{t}^{\cH}  \bD_{0,t}^{\cH} \bB_{mn}^{\cH} \be_{mn}^{\phi} = \trace\{\bB_{mn} \bD_{0,t} \bX_t  \bD_{0,t}^{\cH} \bB_{mn}^{\cH} \be_{mn}^{\phi} (\be_{mn}^{\phi})^{\cT} \} = \trace\{\Phi_{I,mn,t}^{\phi} \bX_t\}$.

\emph{Proof of Lemma}~\ref{SDPformPF}. After standard manipulations,~\eqref{eq:V} can be re-written as
$(P_{L,n,t}^{\phi} + Q_{L,n,t}^{\phi} - y_{C,n}^{\phi}\trace(\bPhi_{V,n}^{\phi} \bX_t))^2 \leq (P_{L,n,t}^{\phi} \eta_{n,t}^{-1})^2$, which is quadratic in $\bX_t$. Then,~\eqref{eq:SDPformPF} is readily obtained by using Schur's complement.

%%%%%%%%%%%%%%%%%%%%%%%%%%%%%%%%%%%%%%%%%%%%%%
\bibliographystyle{IEEEtran}
\bibliography{Biblio_pow_systems}
%%%%%%%%%%%%%%%%%%%%%%%%%%%%%%%%%%%%%%%%%%%%%%

%%%%%%%%%%%%%%%%%%%%%
%\begin{figure}
%\begin{center}
%\includegraphics[width=0.65\textwidth]{F_Scenario}
%\caption{Test scenario.}
%\label{fig:Scenario}
%\end{center}
%\end{figure}
%%%%%%%%%%%%%%%%%%%%%

\vspace{-1.2cm}

\begin{biographynophoto}
{Emiliano Dall'Anese (S'08, M'11)} received the Laurea Triennale (B.Sc degree) and the Laurea Specialistica (M.Sc degree) in Telecommunications Engineering from the University of Padova, Italy, in 2005 and 2007, respectively, and the Ph.D in Information Engineering at the Department of Information Engineering (DEI), University of Padova, Italy, in 2011. From January 2009 to September 2010 he was a visiting scholar at the Department of Electrical and Computer Engineering, University of Minnesota, USA. He is currently a post-doctoral associate at the Department of Electrical and Computer Engineering, University of Minnesota, USA. 

His research interests lie in the areas of signal processing, communications, and networking. Current research focuses on wireless cognitive networks, IP networks,  and power distribution systems.
\end{biographynophoto}

\vspace{-1.2cm}

\begin{biographynophoto}{Georgios B. Giannakis (F'97)}
received his Diploma in Electrical 
Engr. from the Ntl. Tech. Univ. of Athens, Greece, 1981. From 
1982 to 1986 he was with the Univ. of Southern California (USC), 
where he received his MSc. in Electrical Engineering, 1983, MSc. 
in Mathematics, 1986, and Ph.D. in Electrical Engr., 1986. Since 
1999 he has been a professor with the Univ. of Minnesota, where 
he now holds an ADC Chair in Wireless Telecommunications in the 
ECE Department, and serves as director of the Digital Technology
Center. 

His general interests span the areas of communications, networking and
statistical signal processing - subjects on which he has published more
than 300 journal papers, 500 conference papers, 20 book chapters, two
edited books and two research monographs. Current research focuses on
compressive sensing, cognitive radios, network coding, cross-layer designs,
wireless sensors, social and power grid networks. 
\end{biographynophoto}

\vspace{-1.2cm}

\begin{biographynophoto}{Bruce F. Wollenberg (M'67, SM'75, F'89, LF'08)} 
graduated from Rensselaer Polytechnic Institute, Troy, NY,
with a B.S.E.E. in 1964 and an M.Eng. in electric
power engineering in 1966. He then attended
the University of Pennsylvania, Philadelphia, PA,
graduating with a Ph.D. in systems engineering in
1974. He worked for Leeds and Northrup Co. North
Wales, PA from 1966 to 1974, Power Technologies
Inc. Schenectady, NY from 1974 to 1984, and Control
Data Corporation Energy Management System
Division Plymouth, MN from 1984 to 1989. He
took a position as a Professor of Electrical Engineering in the Electrical
and Computer Engineering Department at the University of Minnesota in
September 1989. He is presently the Director of the University of Minnesota
Center for Electric Energy (UMCEE). His main research interests are the
application of mathematical analysis to power system operation and planning
problems. He is the co-author of the textbook \emph{Power Generation Operation
and Control} published by John Wiley $\&$ Sons
\end{biographynophoto}

\end{document}